\documentclass[12pt]{amsart}
\usepackage{epsfig}
\usepackage{graphics}
\numberwithin{equation}{section}
\theoremstyle{plain}
   \newtheorem{thm}{Theorem}[section]
   \newtheorem{cor}[thm]{Corollary}
   
   \newtheorem{lemma}[thm]{Lemma}

\theoremstyle{definition}
   \newtheorem{dfn}[thm]{Definition}
   
\theoremstyle{remark}
   \newtheorem{rem}[thm]{Remark}

\newcommand{\R}{\mathbb R}

\newcommand{\ve}{\varepsilon}

\newcommand{\TB}{\text{\textit{TB}}}
\newcommand{\tb}{\text{\textit{tb}}}

\begin{document}
\title{On the Thurston-Bennequin invariant of graph divide links}
\author{Masaharu Ishikawa}
\footnote[0]{
  This work is supported by the Sumitomo Foundation,
  Grant for Basic Japanese Research Project.
  It is also supported by MEXT,
  Grant-in-Aid for Young Scientists (B) (No. 16740031).
}
\address{Department of Mathematics, Tokyo Institute of Technology,\\
2-12-1, Oh-okayama, Meguro-ku, Tokyo, 152-8551, Japan}
\email{ishikawa@math.titech.ac.jp}

\begin{abstract}
In the present paper we determine the Thurston-Bennequin invariant
of graph divide links,
which include all closed positive braids, all divide links
and certain negative twist knots.
As a corollary of this and a result of
P.~Lisca and A.I.~Stipsicz, we prove that
the $3$-manifold obtained from $S^3$ by Dehn surgery along
a non-trivial graph divide knot $K$ with coefficient $r$
carries positive, tight contact structures for every $r$
except the Thurston-Bennequin invariant of $K$.
\end{abstract}

\maketitle

\section{Introduction}

A contact form on an oriented, closed, smooth $3$-manifold $M$
is an $1$-form $\alpha$ satisfying $\alpha\land d\alpha\ne 0$.
The $2$-plane field $\xi$ on $M$ determined by the kernel of $\alpha$
is called a {\it contact structure}.
A contact structure is called {\it positive}
if the orientation determined by the contact structure
agrees with the given one on $M$.

Let $(M,\xi)$ be a contact manifold, i.e.,
an oriented, closed, smooth $3$-manifold $M$ equipped with
a contact structure $\xi$.
A {\it Legendrian curve} $K$ in $M$ is a simple closed curve in $M$
such that at each point in $K$ the tangent vector to $K$
is included in the $2$-plane of $\xi$.
An {\it overtwisted disk} is a disk in $M$ whose boundary
is a Legendrian curve $K$ and at each point in $K$
the disk is transverse to the $2$-plane of $\xi$.
A contact structure $\xi$ is called {\it tight}
if there are no overtwisted disks in $M$.

It is well-known that $S^3$ admits a unique, positive,
tight contact structure $\xi_0$,
so-called the {\it standard contact structure}.
Let $K$ be a Legendrian link in $(S^3,\xi_0)$, which is an oriented link
consisting of Legendrian curves in $(S^3,\xi_0)$,
and denote by $K^\perp$ an oriented link obtained by pushing-off $K$
in the direction normal to the $2$-plane field $\xi_0$.
The linking number $lk(K,K^\perp)$ of $K$ and $K^\perp$
is called the {\it Thurston-Bennequin invariant} of $K$
and we denote it by $\tb(K)$.
It is known by D.~Bennequin that $\tb(K)$ has an upper bound
when $K$ runs in the same ambient isotopy type~\cite{bennequin:83}.
The maximal number of $\tb(K)$ is called
the {\it maximal Thurston-Bennequin invariant} of $K$
and denoted by $\TB(K)$.

In the present paper we will determine
the maximal Thurston-Bennequin invariant
for a class of links, called {\it graph divide links}.
A {\it graph divide}, introduced by T.~Kawamura in~\cite{kawamura},
is the image of a generic immersion $\varphi:G\to D$
of a disjoint sum $G$ of intervals, circles and finite graphs
into the unit disk $D$.
We call the image of a vertex of $G$ a {\it vertex} of $\varphi(G)$
and a transversal crossing point produced by the immersion
a {\it double point} of $\varphi(G)$.

\begin{figure}[htbp]
   \centerline{\input{tb3.pstex_t}}
   \caption{An example of a graph divide and its knot. This knot is $8_{21}$.
(see \cite[Example 3.4]{kawamura}.)
\label{tb3}}
\end{figure}

Let $P=\varphi(G)\subset D$ be a graph divide.
Suppose that each vertex of $P$ has a sign $+$ or $-$.
We then define an oriented link using $P$ and its signs, see Section~2.
We call it the link of a graph divide $P$ with signs, or shortly
a graph divide link of $P$ (with signs).
It is worth noting that the class of graph divide links includes,
for example,
\begin{itemize}
   \item[(i)] all closed positive braid (see~\cite[Remark 5.4]{kawamura}),
   \item[(ii)] all divide links (see~\cite{acampo:99,acampo:98a}
   and the list in~\cite{gi2}), and
   \item[(iii)] twist knots with $q<0$ twists
                in the notation in~\cite[p.112]{rolfsen}
               (see~\cite{gibson2} and \cite[Example 3.3]{kawamura},
                and cf.~the first footnote in~\cite{ls}). 
\end{itemize}
Note that the links of complex plane curve singularities
are included in both (i) and (ii).

The slice Euler characteristic $\chi_s(K)$
of an oriented link $K$ in $S^3$ is the maximal number of
Euler characteristics of 2-dimensional manifolds in $B^4$,
bounded by $S^3$, whose boundary is $K$.
The integer $g_s(K)=(2-\mu-\chi_s(K))/2$ is called the {\it $4$-genus} of $K$,
where $\mu$ is the number of components of the link $K$.

The main theorem in the present paper is the following:

\begin{thm}\label{thm1}
Let $K$ be an oriented link obtained as
the link of a graph divide $P=\varphi(G)$ with signs $\ve$.
Then we have $\TB(K)=-\chi_s(K)$.
\end{thm}


Before stating a corollary, we introduce a result of 
P.~Lisca and A.I.~Stipsicz.

\begin{thm}[Lisca and Stipsicz~\cite{ls}]\label{thmls}
Let $K$ be a knot in $S^3$ such that $g_s(K)>0$ and $\TB(K)=2g_s(K)-1$. 
Then the oriented $3$-manifold obtained from $S^3$
by Dehn surgery along $K$ with coefficient $r$ carries
positive, tight contact structures for every $r\ne\TB(K)$.
\end{thm}

Since $\TB(K)=-\chi_s(K)=2g_s(K)-1$ for every graph divide knot $K$,
we have the following corollary:

\begin{cor}\label{cor1}
Let $K$ be a non-trivial knot obtained as
the link of a graph divide $P=\varphi(G)$ with signs $\ve$.
Then the oriented $3$-manifold obtained from $S^3$
by Dehn surgery along $K$ with coefficient $r$ carries
positive, tight contact structures for every $r\ne\TB(K)$.
\end{cor}

We will recall the definition of the link of a graph divide in Section~2
and prove Theorem~\ref{thm1} and Corollary~\ref{cor1} in Section~3.

The author would like to thank Tomomi Kawamura for her useful suggestions.

\section{Graph divides}

\begin{figure}[htbp]
   \centerline{\input{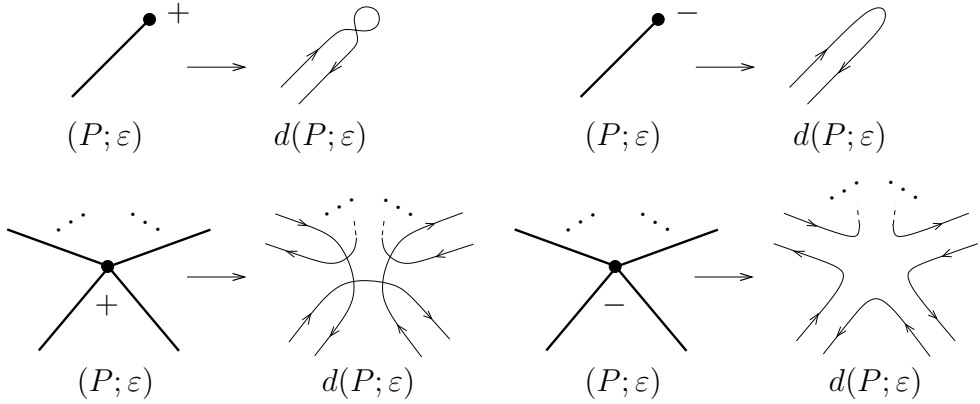}}
   \caption{Doubling method around vertices of graph divides.\label{tb2}}
\end{figure}

Let $P=\varphi(G)\subset D$ be a graph divide.
Suppose that each vertex of $P$ has a sign $+$ or $-$.
We denote this distribution of signs by $\ve$.
To each pair $(P,\ve)$ of a graph divide $P$
and a distribution $\ve$ of signs,
we define oriented, immersed circles in $D$ according to
the following {\it doubling method}:
\begin{itemize}
   \item[(i)] outside a small neighborhood of the vertices of $P$
we draw two curves parallel to $P$;
   \item[(ii)] inside the small neighborhood of each vertex of $P$
we draw the curves as shown in Figure~\ref{tb2}
according to the sign at the vertex;
   \item[(iii)] assign an orientation to the curves obtained
so that it is ``keeping left'' almost everywhere.
\end{itemize}
We denote the obtained oriented, immersed circles by $d(P;\ve)$.

An {\it oriented divide} is the image of a generic immersion of
oriented circles into $D$~\cite{gi}.
So, $d(P;\ve)$ is an oriented divide.
Let $T(\R^2)$ be the tangent bundle to $\R^2$,
on which the unit disk $D$ lies,
and $ST(\R^2)$ the unit sphere in $T(\R^2)$.
We denote by $x=(x_1,x_2)$ the coordinates of $\R^2$
and by $u=(u_1,u_2)$ the coordinates of tangent space to $\R^2$.
The {\it link} of an oriented divide $\vec P$ is the set 
\[
   L(\vec P)=\{(x,u)\in ST(\R^2)\mid x\in \vec P,\,u\in T_x(\vec P)\}
   \subset S^3,
\]
where $x\in \vec P$ means that $x$ lies on the curve of $\vec P$
and $T_x(\vec P)$ is the set of tangent vectors to $\vec P$
whose direction is consistent with the orientation assigned to $\vec P$.

\begin{dfn}
The {\it link} $L(P;\ve)$ of a graph divide $P$ with signs $\ve$
is the set given by
\[
   L(P;\ve)=\{(x,u)\in ST(\R^2)\mid x\in d(P;\ve),\,u\in T_x(d(P;\ve))\}
   \subset S^3,
\]
i.e., it is the link of the oriented divide $d(P;\ve)$.
\end{dfn}

\begin{rem}
We here note a brief history. $P=\varphi(G)$ is called a {\it divide}
if $G$ is a disjoint sum of intervals and circles
and the image of the endpoints of intervals lies on the boundary of $D$.
The divide is introduced by N.~A'Campo
as a generalization of real morsified curves of
complex plane curve singularities.
In particular, their links are fibered and the fibrations correspond
to the Milnor fibrations of the complex
singularities.
He also determined the unknotting numbers and the slice Euler
characteristics of all divide knots~\cite{acampo:99,acampo:98a}.
Next W.~Gibson and the author introduced {\it free divides}
by forgetting the assumption about the endpoints of intervals~\cite{gi2}.
After that, W.~Gibson introduced {\it tree divides}
as the image of trees and
determined their unknotting numbers, slice Euler characteristics, and
Thurston-Bennequin invariant~\cite{gibson2}.
Finally, T.~Kawamura introduced graph divides, proved
their quasipositivity and determined their
slice Euler characteristics~\cite{kawamura}.
\end{rem}

\section{Proofs}

Let $P=\varphi(G)$ be a graph divide with signs $\ve$
and consider the $\pi/2$-rotation map $R_{\pi/2}:ST(\R^2)\to ST(\R^2)$
given by $(x_1,x_2,u_1,u_2)\mapsto (x_1,x_2,u_2,-u_1)$.
This map corresponds to the $\pi/2$-rotation of tangent vectors
to the unit disk $D$
and the image of the $\pi/2$-rotation of tangent vectors to $d(P;\ve)$
constitutes co-oriented, immersed circles $w$ in $D$, i.e.,
a wave front on $D$.
Since $R_{\pi/2}(L(P;\ve))$ is the preimage of the wave front $w$,
it is Legendrian with respect to the standard contact structure in $S^3$.
In what follows, we regard $L(P;\ve)$ as a Legendrian link
up to the $\pi/2$-rotation map $R_{\pi/2}$.

Let $\delta(P)$ denote the number of double points of $P=\varphi(G)$
and $\chi(G)$ the Euler characteristic of $G$.
We will prove the following:

\begin{thm}\label{thm1a}
For a graph divide $P=\varphi(G)$ with signs $\ve$,
we have $\tb(L(P;\ve))=2\delta(P)-\chi(G)$.
\end{thm}

Before proving Theorem~\ref{thm1a} we give a proof of the main theorem.

\vspace{3mm}

\noindent
{\it Proof of Theorem~\ref{thm1}.}\,
Let $K$ denote the link $L(P;\ve)$ without information about Legendrian.
It is known in~\cite[Proposition 6.2]{kawamura} that
$\chi_s(K)=\chi(G)-2\delta(P)$.
Hence we have $\tb(L(P;\ve))=-\chi_s(K)$.
By the slice Bennequin inequality~\cite{rudolph},
we have $\tb(K^{Leg})\leq -\chi_s(K)$
for any Legendrian link $K^{Leg}$ ambient isotopic to $K$
and hence $\TB(K)=\tb(L(P;\ve))=-\chi_s(K)$.
\qed

\vspace{3mm}

In order to prove Theorem~\ref{thm1a}
we introduce a modified tangle presentation of a graph divide.
First prepare seven kinds of tangles as shown in Figure~\ref{tb1}.
We name them as ``a left endpoint tangle'',
``a right endpoint tangle'', and so on.
The left or right endpoint tangle has two cases depending on
the sign $+$ or $-$ at the endpoint vertex.
A {\it tangle product} is a product of these tangles with well-defined
connections.

\begin{figure}[htbp]
   \centerline{\input{tb1.pstex_t}}
   \caption{Tangles.\label{tb1}}
\end{figure}

\begin{dfn}
A {\it modified tangle presentation} of a graph divide $P=\varphi(G)$
with signs $\ve$ is a tangle product whose immersed graph has
the same configuration as the graph divide $P$ with signs $\ve$.
\end{dfn}

It is easy to see that any graph divide $P$ with signs $\ve$
has a modified tangle presentation.

Let $e_1$ (resp. $e_2$) be the number of left (resp. right) endpoint tangles
and $f_1$ (resp. $f_2$) the number of left (resp. right) fold tangles.
We denote by $v_1,\cdots,v_m$ the vertices in the ($+$) branch tangles
and by $b(v_i)$ the number of edges connected to the vertex $v_i$.

\begin{lemma}\label{lemma1}
$-\chi(G)=\sum_{i=1}^m(b(v_i)-2)-e_2+f_1-f_2$.
\end{lemma}

\begin{proof}
By checking $\chi(G)$ from the right to the left
in the modified tangle presentation,
we have $\chi(G)=e_2+f_2-\sum_{i=1}^m(b(v_i)-2)-f_1$.
\end{proof}

\begin{figure}[htbp]
   \centerline{\input{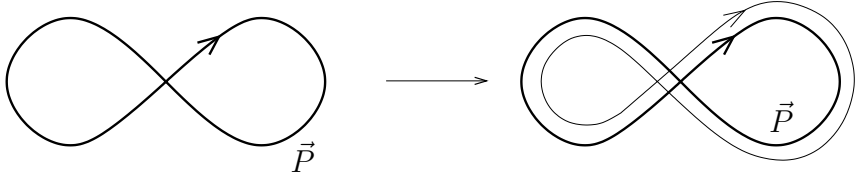}}
   \caption{A parallel shift of an oriented divide $\vec P$.\label{tb4}}
\end{figure}

Next we calculate $\tb(L(P;\ve))=lk(L(P;\ve),L^\perp(P;\ve))$.
The link $R_{\pi/2}(L^\perp(P;\ve))$ is the Legendrian
link of a wave front $w^\perp$ obtained by pushing-off the wave front $w$ 
in the direction of the assigned co-orientation.
Since $R_{\pi/2}$ is an isomorphism, the linking of two links does not change.
So, $lk(L(P;\ve),L^\perp(P;\ve))$ is equal to
the linking number of the Legendrian link of $w$ and that of $w^\perp$.
By the $-\pi/2$-rotation map $R^{-1}_{\pi/2}$, we conclude that
$lk(L(P;\ve),L^\perp(P;\ve))$ is equal to the linking number
of $L(P;\ve)$ and the link of an oriented divide obtained
by parallel-shifting $d(P;\ve)$ in one of the possible directions,
see Figure~\ref{tb4}.

A method of obtaining a link diagram of the link of an oriented divide
has been made in~\cite{gi}.
Let $\vec P$ be an oriented divide in $D$.
Each point $(x,u)$ in $L(\vec P)$ corresponds to a point $x\in D$
and a vector $u\in T_x(\vec P)$.
Fix a parametrization $x(t)$ of $\vec P$ with parameter
$t\in \sqcup_{i=1}^\mu S^1$, where $\mu$ is the number of
components of $\vec P$.
This induces a parametrization $(x(t),u(t))$ of the link $L(\vec P)$.
By setting $\vec P$ in general position,
we assume that $y=\arg u(t)$ is transversal to $y\equiv -\pi/2 \mod 2\pi$
and if $x(t_1)=x(t_2)$ for $t_1\ne t_2$ then
$\arg u(t_i)\not\equiv -\pi/2 \mod 2\pi$ for $i=1,2$.

We regard the curve of oriented divide $\vec P$
as the strands of the link $L(\vec P)$
by applying the following rules:

\begin{dfn}[over/under rule]
Suppose that $x(t_1)=x(t_2)$ for $t_1\ne t_2$.
Let $\theta_i$,\,$i=1,2$,\,be real numbers such that
$-\pi/2 < \theta_i < 3\pi/2$ and $\theta_i\equiv \arg u(t_i)\mod 2\pi$.
If $\theta_1 < \theta_2$ then the strand corresponding to $t_1$ passes
over the strand corresponding to $t_2$.
\end{dfn}

We are drawing the link diagram on $\R^2\supset D$ with
coordinates $(x_1,x_2)$.
Let $x_3$ denote the third coordinate for describing
over/under at crossing points of the link diagram.

\begin{dfn}[winding rule]
If $\arg u(t_0)\equiv -\pi/2 \mod 2\pi$ and $\arg u(t)$ is locally increasing
(resp. decreasing) at $t=t_0$ then we cut the strand at that point, bring it
to $x_3=-\infty$ (resp. $x_3=\infty$)
and take it back from $x_3=\infty$ (resp. $x_3=-\infty$).
\end{dfn}

We fix the coordinates $(x_1,x_2)\in\R^2$ in such a way that
$x_1$ is in the right direction in Figure~\ref{tb1}
and $x_2$ is in the top.
In order to calculate $lk(L(P;\ve),L^\perp(P;\ve))$,
we draw the link diagram of the oriented divide $d(P;\ve)$
and its parallel-shift, and check the algebraic number of crossings
of these links.

\begin{lemma}\label{lemma2}
The contribution to $lk(L(P;\ve),L^\perp(P;\ve))$ in each tangle
is as follows:
\begin{itemize}
   \item[(i)] $0$ for each left endpoint tangle;
   \item[(ii)] $-1$ for each right endpoint tangle;
   \item[(iii)] $+1$ for each left fold tangle;
   \item[(iv)] $-1$ for each right fold tangle;
   \item[(v)] $b(v_i)-2$ for each {\rm ($+$)} branch tangle with vertex $v_i$;
   \item[(vi)] $0$ for each {\rm ($-$)} branch tangle;
   \item[(vii)] $+2$ for each double point tangle.
\end{itemize}
\end{lemma}

\begin{figure}[htbp]
   \centerline{\input{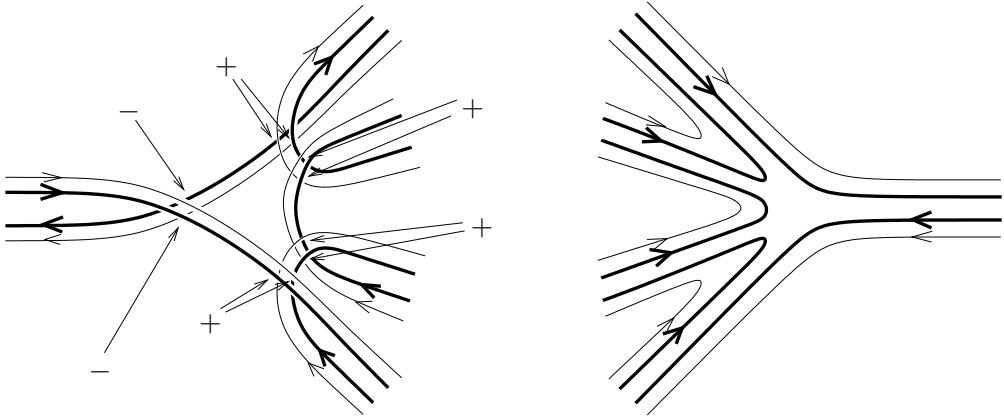}}
   \caption{Link diagrams for ($+$) and ($-$) branch tangles.\label{tb8}}
\end{figure}

\begin{proof}
The cases (i) $\sim$ (iv) and (vii) have been proved
by W.~Gibson in~\cite{gibson}.
We here write their proofs again for completeness of the proof.

Let $\vec P$ denote the oriented divide consisting of
$d(P;\ve)$ and its parallel-shift.
First we consider the case of a left endpoint tangle with $+$ sign.
The corresponding link diagram is described in Figure~\ref{tb5}.
In the left endpoint tangle, two positive crossings and
two negative crossings appear and hence their algebraic sum is $0$.
Outside the tangle the strand drawn by the winding rule may produce
crossings with other strands. However since $\vec P$ is
obtained by the doubling method
the algebraic sum of these signs becomes always $0$,
see the crossings around the four parallel strands in Figure~\ref{tb5}.
Hence the contribution of this tangle to $lk(L(P;\ve),L^\perp(P;\ve))$ is $0$.

\begin{figure}[htbp]
   \centerline{\input{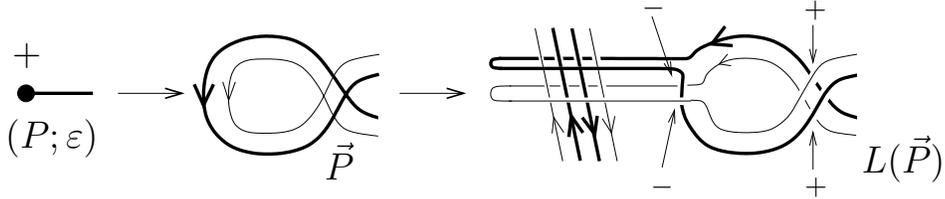}}
   \caption{A link diagram for a left endpoint tangle with sign $+$.\label{tb5}}
\end{figure}

If a left endpoint tangle has $-$ sign then the oriented divide $\vec P$
has no double points and no point $x(t)$ at
which $\arg u(t)\equiv -\pi/2 \mod 2\pi$.
Hence the corresponding link diagram has no crossings.
This completes the proof of case (i).

The case of a right endpoint tangle with sign $+$
is described in Figure~\ref{tb6}. Hence the contribution is $-1$.
In the case of $-$ sign,
the winding rule produces two negative crossings
and the contribution becomes also $-1$.
This completes the proof of case (ii).

\begin{figure}[htbp]
   \centerline{\input{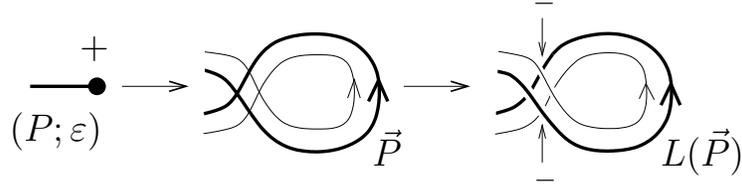}}
   \caption{A link diagram for a right endpoint tangle with sign $+$.\label{tb6}}
\end{figure}

The cases of left/right fold tangles and a double point tangle
are described in Figure~\ref{tb7}.
The windings outside the fold tangles do not contribute to the calculation
of the linking number because of the same reason as we mentioned
in case (i).
The contributions are $+1$, $-1$ and $+2$ respectively
and these complete the proofs of cases (iii), (iv) and (vii).

\begin{figure}[htbp]
   \centerline{\input{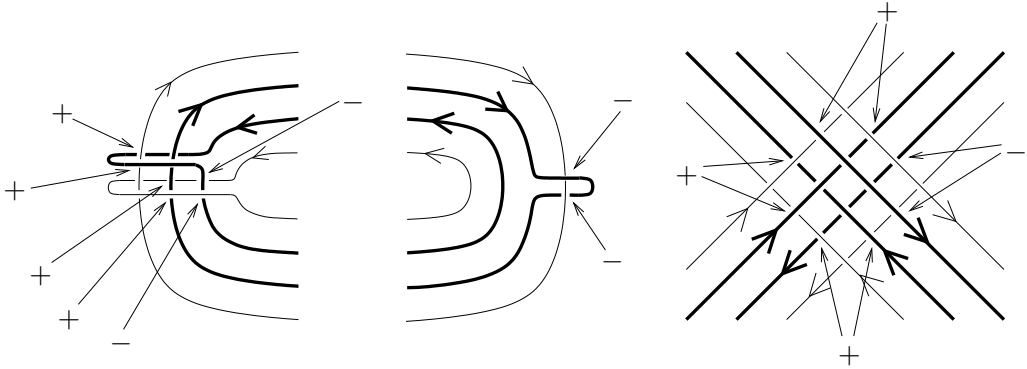}}
   \caption{Link diagrams for left/right fold tangles and a double point tangle.\label{tb7}}
\end{figure}

The cases of branch tangles are described in Figure~\ref{tb8}.
The figures are in the case where the numbers of branches
at the vertices are $5$.
Thus the contribution is $b(v_i)-2$ if it is a ($+$) branch tangle
and $0$ if it is ($-$).
This completes the proof.
\end{proof}

\noindent
{\it Proof of Theorem~\ref{thm1a}}.\,
From Lemma~\ref{lemma2} we have
$lk(L(P;\ve),L^\perp(P;\ve))=\sum_{i=1}^m(b(v_i)-2)-e_2+f_1-f_2+2\delta(P)$
and this coincides with $2\delta(P)-\chi(G)$ by Lemma~\ref{lemma1}.
\qed

\vspace{3mm}

\noindent
{\it Proof of Corollary~\ref{cor1}}.\,
Since we already have $\TB(K)=-\chi_s(K)$,
it is enough to show that $g_s(K)>0$.
Assume $g_s(K)=0$. Then, by~\cite[Proposition~6.2]{kawamura},
we have $\delta(P)=0$ and $\chi(G)=1$.
But in this case $P$ is an embedded tree and therefore
the link $L(P;\ve)$ is a trivial knot (see~\cite[Lemma~6.8]{kawamura}).
This contradicts the assumption.
\qed

\vspace{3mm}

\begin{rem}
A {\it positive link} is an oriented link which has a positive diagram,
i.e. a diagram with only positive crossings.
Using the results in~\cite{tanaka} and~\cite{nakamura},
we can show the equality $TB(K)=-\chi_s(K)$ for every positive link $K$.
Hence the oriented $3$-manifold obtained from $S^3$
by Dehn surgery along a positive knot $K$ with coefficient $r$
carries positive, tight contact structures for every $r\ne\TB(K)$.
\end{rem}

\begin{rem}
The knot $8_{20}$ is slice and quasipositive
(cf.~\cite[Proof of Theorem 6.7]{kawamura})
and satisfies $\TB(8_{20})=-2$ (see~\cite{ng}).
Hence the equality $\TB(K)=2g_s(K)-1$ does not hold.
This means that the equality is a property of graph divide knots and
positive knots, but is not a property of quasipositive knots.
We don't know if it is a property of strongly quasipositive knots or not.
We remark that the graph divide knot $8_{21}$ shown in Figure~\ref{tb3}
is not strongly quasipositive, see~\cite[Remark 6.9]{kawamura}.
\end{rem}

\end{document}